\documentclass[11pt, twoside]{article}
\usepackage{latexsym}
\usepackage{amsmath}
\usepackage{amssymb}
\usepackage[all]{xy}
\usepackage{amsfonts}
\usepackage{verbatim}
\usepackage{amsthm}
\usepackage{mathrsfs}
\usepackage{epsfig}
\usepackage{xy}
\usepackage{array}
\usepackage{stmaryrd}
\usepackage{graphicx,color}
\usepackage{xcolor}
\usepackage{tikz}
\usetikzlibrary{arrows,calc}
\usepackage{etex}
\usepackage{mathdots}
\usepackage{float}
\usepackage{graphics}
\usepackage{pdflscape}

\usepackage{anysize,hyperref}
\input xypic
\xyoption{all}

\usepackage[perpage,symbol]{footmisc}
\usepackage{setspace}
\def\D{\mathcal{D}}
\def\C{\mathcal{C}}

\def\dr{\ar@{->}[r]}

\def\X{\mathscr{X}}

\newcommand{\add}{\mathsf{add}\hspace{.01in}}
\def\mod{\mathsf{mod}\hspace{.01in}}
\def\proj{\mathsf{proj}\hspace{.01in}}
\def\End{\mathsf{End}}

\def\Hom{\mbox{Hom}}

\usepackage{tikz} \usetikzlibrary{calc} \usepackage{tikz-cd}

\begin{document}
\baselineskip=15pt
\title{\Large{\bf Modules of infinite projective dimension$^\bigstar$\footnotetext{\hspace{-1em}$^\bigstar$This work was supported by the National Natural Science Foundation of China (Grants No. 11901190 and 11671221), and the Hunan Provincial Natural Science Foundation of China (Grants No. 2018JJ3205), and the Scientific Research Fund of Hunan Provincial Education Department (Grant No. 19B239).}}}
\medskip
\author{Panyue Zhou and Xingjia Zhou}

\date{}

\maketitle
\def\blue{\color{blue}}
\def\red{\color{red}}

\newtheorem{theorem}{Theorem}[section]
\newtheorem{lemma}[theorem]{Lemma}
\newtheorem{corollary}[theorem]{Corollary}
\newtheorem{proposition}[theorem]{Proposition}
\newtheorem{conjecture}{Conjecture}
\theoremstyle{definition}
\newtheorem{definition}[theorem]{Definition}
\newtheorem{question}[theorem]{Question}
\newtheorem{remark}[theorem]{Remark}
\newtheorem{remark*}[]{Remark}
\newtheorem{example}[theorem]{Example}
\newtheorem{example*}[]{Example}
\newtheorem{condition}[theorem]{Condition}
\newtheorem{condition*}[]{Condition}
\newtheorem{construction}[theorem]{Construction}
\newtheorem{construction*}[]{Construction}

\newtheorem{assumption}[theorem]{Assumption}
\newtheorem{assumption*}[]{Assumption}

\baselineskip=17pt
\parindent=0.5cm
\vspace{-6mm}

\begin{abstract}
\baselineskip=16pt
We characterize the modules of infinite projective dimension
over the endomorphism algebras of Opperman-Thomas cluster tilting objects $X$ in $(n+2)$-angulated categories $(\C,\Sigma^n,\Theta)$.
For an indecomposable object $M$ of $\C$,
 we define in this article the ideal $I_M$ of $\End_{\C}(\Sigma^nX)$ given by
all endomorphisms that factor through $\add M$, and show that the $\End_{\C}(X)$-module $\Hom_{\C}(X,M)$
has infinite projective dimension precisely when $I_M$ is non-zero. As an application,  we generalize
a recent result by Beaudet-Br\"{u}stle-Todorov for cluster-tilted algebras.\\[0.2cm]
\end{abstract}

\pagestyle{myheadings}
\markboth{\rightline {\scriptsize   Panyue Zhou and Xingjia Zhou}}
         {\leftline{\scriptsize Modules of infinite projective dimension}}

\section{Introduction}
Cluster-tilting theory sprouted from the categorification of Fomin-Zelevinsky's cluster algebras. It is used to construct abelian categories from some
triangulated categories. By Buan-Marsh-Reiten \cite[Theorem 2.2]{BMR} in the case of cluster categories, by Keller-Reiten
\cite[Proposition 2.1]{KR} in the $2$-Calabi-Yau case, then by Koenig-Zhu \cite[Theorem 3.3]{KZ} and Iyama-Yoshino \cite[Corollary 6.5]{IY} in the general case, one can pass from triangulated categories to abelian categories by factoring out cluster tilting subcategories.

Recall that the notion of cluster-tilting subcategory which is due to Iyama \cite{I}.
Let $\C$ be a triangulated category with shift functor $\Sigma$. A subcategory $\X$ of $\C$ is called  \emph{cluster tilting} if
it satisfies the following conditions:

(1) $\X$ is contravariantly finite and covariantly finite in $\C$.

(2) $X\in\X$ if and only if $\Hom_{\C}(X,\Sigma\X)=0$, i.e. $\Hom_{\C}(X,\Sigma M)=0,$ for any $M\in \X$.

(3) $X\in\X$ if and only if $\Hom_{\C}(\X,\Sigma X)=0$,  i.e. $\Hom_{\C}(M,\Sigma X)=0,$ for any $M\in \X$.

An object $X$ is called \emph{cluster tilting}, if $\add X$ is cluster tilting,  where $\add X$ is the subcategory of $\C$ consisting of direct summands of direct sum of finitely many copies of $X$.

In fact, Koenig and Zhu \cite[Lemma 3.2]{KZ} show that $\X$ is cluster tilting if and only if $\Hom_{\C}(\X,\Sigma\X)=0$, i.e. $\Hom_{\C}(X,\Sigma Y)=0$, for any $X, Y\in \X$, and for any object $C\in\C$, there exists a triangle
$X_0\to X_1\to C\to \Sigma X_0$
where $X_0,X_1\in\X$.

Let $\C$ be a triangulated category with a shift functor $\Sigma$, and $X\in\C$ be a cluster tilting object. For any object $M\in\C$, we denote by
$I_M$ the ideal of $\Gamma:=\End_{\C}(\Sigma X)$ ($\Gamma$ is called a cluster-tilted algebra) given by all endomorphisms that factor through $\add M$ and call it
\emph{factorization ideal} of $M$.
Beaudet-Br\"{u}stle-Todorov showed the following.
\begin{theorem}\emph{\cite[Theorem 1.1]{BBT}}
Let $\C$ be a triangulated category with a shift
functor $\Sigma$ and $X\in\C$ be a cluster tilting object. Let $M$ be an indecomposable object in $\C$ which does not belong to $\add \Sigma X$. Then the $\End_{\C}(X)$-module ${\rm Hom}_{\C}(X,M)$ is of infinite
projective dimension if and only if the factorization ideal $I_M$ is non-zero.
\end{theorem}

Later, Liu and Xiu \cite[Theorem 3.1]{LX}  generalized the work of Beaudet-Br\"{u}stle-Todorov for cluster-tilted algebras to the endomorphism algebras of
maximal rigid objects in $2$-Calabi-Yau triangulated categories.

Recently, Geiss, Keller and Oppermann introduced in \cite{GKO} a new type of categories, called $(n+2)$-angulated categories, which generalize triangulated categories: the classical triangulated categories are the special case $n=1$. These categories appear for instance
when considering certain $n$-cluster tilting subcategories of triangulated categories.

The notion of cluster tilting objects can be generalised to $(n+2)$-angulated categories, which is due to Oppermann and Thomas \cite[Definition 5.3]{OT}.

\begin{definition}\cite[Definition 5.3]{OT}
Let $\C$ be an $(n+2)$-angulated category with an  $n$-suspension
functor $\Sigma^n$. An object $X\in\C$ is called \emph{Opperman-Thomas cluster tilting}  if

(1) $\Hom_{\C}( X,\Sigma^n X)=0$.

(2)  For any object $C\in\C$, there exists an $(n+2)$-angle
$$X_0\xrightarrow{~~}X_1\xrightarrow{~~}\cdots\xrightarrow{~~}X_n\xrightarrow{~~}C\xrightarrow{~~}\Sigma^n X_0$$
where $X_0, X_1,\cdots,X_{n}\in\add X$.
\end{definition}

Let $\C$ be an $(n+2)$-angulated category with an $n$-suspension
functor $\Sigma^n$, and $X\in\C$ be an Opperman-Thomas cluster tilting object. For any object $M\in\C$, we denote by
$I_M$ the ideal of $\End_{\C}(\Sigma^n X)$ given by all endomorphisms that factor through $M$ and call it
\emph{factorization ideal} of $M$.
Our main result is the following, which is a
generalization of Beaudet-Br\"{u}stle-Todorov result.
\begin{theorem}\emph{(See Theorem \ref{main} for details)}
Let $\C$ be an $(n+2)$-angulated category with an $n$-suspension
functor $\Sigma^n$ and $X\in\C$ be an Opperman-Thomas cluster tilting object. Let $M$ be an indecomposable object in $\C$ which does not belong to $\add \Sigma^nX$. Then the $\End_{\C}(X)$-module ${\rm Hom}_{\C}(X,M)$ is of infinite
projective dimension if and only if the factorization ideal $I_M$ is non-zero.
\end{theorem}

This article is organised as follows: In Section 2, we review some elementary definitions and facts
of $(n+2)$-angulated categories.
In Section 3,  we prove our main result.

\section{Preliminaries}
In this section, we briefly recall the definition and basic properties of $(n+2)$-angulated categories from \cite{GKO}.
Let $\C$ be an additive category with an automorphism $\Sigma^n:\C\rightarrow\C$, and $n$ an integer greater than or equal to one.

An $(n+2)$-$\Sigma^n$-$sequence$ in $\C$ is a sequence of objects and morphisms
$$A_0\xrightarrow{f_0}A_1\xrightarrow{f_1}A_2\xrightarrow{f_2}\cdots\xrightarrow{f_{n-1}}A_n\xrightarrow{f_n}A_{n+1}\xrightarrow{f_{n+1}}\Sigma^n A_0.$$
Its {\em left rotation} is the $(n+2)$-$\Sigma^n$-sequence
$$A_1\xrightarrow{f_1}A_2\xrightarrow{f_2}A_3\xrightarrow{f_3}\cdots\xrightarrow{f_{n}}A_{n+1}\xrightarrow{f_{n+1}}\Sigma^n A_0\xrightarrow{(-1)^{n}\Sigma^n f_0}\Sigma^n A_1.$$
A \emph{morphism} of $(n+2)$-$\Sigma^n$-sequences is  a sequence of morphisms $\varphi=(\varphi_0,\varphi_1,\cdots,\varphi_{n+1})$ such that the following diagram commutes
$$\xymatrix{
A_0 \ar[r]^{f_0}\ar[d]^{\varphi_0} & A_1 \ar[r]^{f_1}\ar[d]^{\varphi_1} & A_2 \ar[r]^{f_2}\ar[d]^{\varphi_2} & \cdots \ar[r]^{f_{n}}& A_{n+1} \ar[r]^{f_{n+1}}\ar[d]^{\varphi_{n+1}} & \Sigma^n A_0 \ar[d]^{\Sigma^n \varphi_0}\\
B_0 \ar[r]^{g_0} & B_1 \ar[r]^{g_1} & B_2 \ar[r]^{g_2} & \cdots \ar[r]^{g_{n}}& B_{n+1} \ar[r]^{g_{n+1}}& \Sigma^n B_0
}$$
where each row is an $(n+2)$-$\Sigma^n$-sequence. It is an {\em isomorphism} if $\varphi_0, \varphi_1, \varphi_2, \cdots, \varphi_{n+1}$ are all isomorphisms in $\C$.

\begin{definition}\cite[Definition 2.1]{GKO}
An $n$-\emph{angulated category} is a triple $(\C, \Sigma^n, \Theta)$, where $\C$ is an additive category, $\Sigma^n$ is an automorphism of $\C$ ($\Sigma^n$ is called the $n$-suspension functor), and $\Theta$ is a class of $(n+2)$-$\Sigma^n$-sequences (whose elements are called $(n+2)$-angles), which satisfies the following axioms:
\begin{itemize}
\item[\textbf{(N1)}]
\begin{itemize}
\item[(a)] The class $\Theta$ is closed under isomorphisms, direct sums and direct summands.

\item[(b)] For each object $A\in\C$ the trivial sequence
$$ A\xrightarrow{1_A}A\rightarrow 0\rightarrow0\rightarrow\cdots\rightarrow 0\rightarrow \Sigma^nA$$
belongs to $\Theta$.

\item[(c)] Each morphism $f_0\colon A_0\rightarrow A_1$ in $\C$ can be extended to an $(n+2)$-angle: $$A_0\xrightarrow{f_0}A_1\xrightarrow{f_1}A_2\xrightarrow{f_2}\cdots\xrightarrow{f_{n-1}}A_n\xrightarrow{f_n}A_{n+1}\xrightarrow{f_{n+1}}\Sigma^n A_0.$$
\end{itemize}
\item[\textbf{(N2)}] An $(n+2)$-$\Sigma^n$-sequence belongs to $\Theta$ if and only if its left rotation belongs to $\Theta$.

\item[\textbf{(N3)}] For each solid commutative diagram
$$\xymatrix{
A_0 \ar[r]^{f_0}\ar[d]^{\varphi_0} & A_1 \ar[r]^{f_1}\ar[d]^{\varphi_1} & A_2 \ar[r]^{f_2}\ar@{-->}[d]^{\varphi_2} & \cdots \ar[r]^{f_{n}}& A_{n+1} \ar[r]^{f_{n+1}}\ar@{-->}[d]^{\varphi_{n+1}} & \Sigma^n A_0 \ar[d]^{\Sigma^n \varphi_0}\\
B_0 \ar[r]^{g_0} & B_1 \ar[r]^{g_1} & B_2 \ar[r]^{g_2} & \cdots \ar[r]^{g_{n}}& B_{n+1} \ar[r]^{g_{n+1}}& \Sigma^n B_0
}$$ with rows in $\Theta$, the dotted morphisms exist and give a morphism of  $(n+2)$-angles

\item[\textbf{(N4)}] In the situation of (N3), the morphisms $\varphi_2,\varphi_3,\cdots,\varphi_{n+1}$ can be chosen such that the mapping cone
$$A_1\oplus B_0\xrightarrow{\left(\begin{smallmatrix}
                                        -f_1&0\\
                                        \varphi_1&g_0
                                       \end{smallmatrix}
                                     \right)}
A_2\oplus B_1\xrightarrow{\left(\begin{smallmatrix}
                                        -f_2&0\\
                                        \varphi_2&g_1
                                       \end{smallmatrix}
                                     \right)}\cdots\xrightarrow{\left(\begin{smallmatrix}
                                        -f_{n+1}&0\\
                                        \varphi_{n+1}&g_n
                                       \end{smallmatrix}
                                     \right)} \Sigma^n A_0\oplus B_{n+1}\xrightarrow{\left(\begin{smallmatrix}
                                        -\Sigma^n f_0&0\\
                                        \Sigma^n\varphi_1&g_{n+1}
                                       \end{smallmatrix}
                                     \right)}\Sigma^nA_1\oplus\Sigma^n B_0$$
belongs to $\Theta$.
   \end{itemize}
\end{definition}

Now we give an example of $(n+2)$-angulated categories.

\begin{example}
We recall the standard construction of an $(n+2)$-angulated category given by Geiss-Keller-Oppermann \cite[Theorem 1]{GKO}.
Let $\C$ be a triangulated category and $\mathcal{T}$ an $n$-cluster tilting subcategory which is closed under $\Sigma^{n}$, where $\Sigma$ is the shift functor of $\C$. Then $(\mathcal{T},\Sigma^{n},\Theta)$ is an $(n+2)$-angulated category, where $\Theta$ is the class of all sequences
$$A_0\xrightarrow{f_0}A_1\xrightarrow{f_1}A_2\xrightarrow{f_2}\cdots\xrightarrow{f_{n-1}}A_n\xrightarrow{f_n}A_{n+1}\xrightarrow{f_{n+1}}\Sigma^{n} A_0$$
such that there exists a diagram
$$\xymatrixcolsep{0.3pc}
 \xymatrix{& A_1 \ar[dr]\ar[rr]^{f_1}  &  & A_2  \ar[dr]  & & \cdots  & & A_{n} \ar[dr]^{f_{n}}      \\
A_0 \ar[ur]^{f_0} & \mid & \ar[ll]  A_{1.5}\ar[ur] & \mid &  \ar[ll]  A_{2.5} & \cdots & A_{n-1.5}\ar[ur] & \mid & \ar[ll] A_{n+1}   }$$
with $A_i\in\mathcal{T}$ for all $i\in\mathbb{Z}$, such that all oriented triangles are triangles in $\C$, all non-oriented triangles commute, and $f_{n+1}$ is the composition along the lower edge of the diagram.
\end{example}

From now on to the end of this section, we assume that $\C$ is an $(n+2)$-angulated category with an $n$-suspension
functor $\Sigma^n$, $X\in\C$ is an Opperman-Thomas cluster tilting object and $\Gamma:=\End_{ \C}( X)$ is the endomorphism algebra of $X$.  Recall the following result.
\begin{theorem}\emph{\cite[Theorem 0.5]{JJ1} and \cite[Theorem 3.8]{ZZ}}\label{y1}
Consider the essential image $\D$ of the functor ${\rm Hom}_\C(X,-) : \C \rightarrow\mod\Gamma$.  Then $\D$ is an $n$-cluster tilting subcategory of $\mod\Gamma$ where $\mod\Gamma$ is the category
of finite dimensional right $\Gamma$-modules.  There exists a commutative diagram, as shown below, where the vertical arrow is the quotient functor and the diagonal arrow is an equivalence of categories:
$$\begin{tikzpicture}[scale=2.5]
\node (T) at (0,1) {\(\C\)};
\node (Tadd) at (0,0) {\(\C/\add\Sigma^n X\).};
\node (D) at (1,1) {\(\D\)};
\draw[->>] (T)--node[left]{} (Tadd);
\draw[->] (T)--node[above]{\(\scriptstyle{{\rm Hom}_\C(X,-)}\)} (D);
\draw[->] (Tadd)-- node[above, sloped] {\(\sim\)} (D);
\end{tikzpicture}$$
Moreover,  $\Gamma$ is an $n$-Gorenstein algebra, that is, each injective module has projective dimension
$\leq n$, and each projective module has injective dimension $\leq n$.
\end{theorem}

\begin{remark}\cite[Lemma 2.1]{JJ1}\label{y2}
The classic $\add$-$\proj$-correspondence holds, as the functor ${\rm Hom}_\C(X,-)$ restricts to an equivalence $\add X\longrightarrow\proj\Gamma$ where
$\proj\Gamma$ is the category of projective finite dimensional right $\Gamma$-modules.
\end{remark}

\begin{corollary}\label{cor}
Each $\Gamma$-module is either of infinite projective dimension or of projective  dimension at most $n$.
\end{corollary}

\proof Suppose that $M$ is a $\Gamma$-module of finite projective dimension and  $M$ is of projective dimension $m$. Then there exists an exact sequence:
\[
  \begin{tikzcd}[column sep=small, row sep=small]
    0\ar{rr}&&\Omega^mM\ar[tail]{rr}{}\drar[equals]&&P_m\ar{rr}{}&&\cdots\ar{rr}{}\drar[two
    heads,swap]{}&&P_0\ar[two heads]{rr}{}\drar[two heads,swap]{}&&M\ar{rr}&&0\\
    &&&\Omega^mM\urar[tail,swap]{}&&&&\Omega M\urar[tail,swap]{}&&M\urar[equals]
  \end{tikzcd}
  \]
where $P_0,P_1,\cdots,P_m$ are projective modules. Since $M$ is of projective dimension $m$, we have that 
$\Omega^mM$ is projective module. By Theorem \ref{y1}, we know that  $\Gamma$ is an $n$-Gorenstein algebra.
It follows that $\Omega^mM$ and $P_m$ are of injective dimension at most $n$.
Thus $\Omega^{m-1}M$ is of injective dimension at most $n$.
By induction on $t\geq 0$, we find that $\Omega^{m-t}M$ is of injective dimension
at most $n$. In particular, $M$ is of injective dimension at most $n$.
Dually, one shows that if $N$ is a $\Gamma$-module of finite injective dimension, then
$N$ is of projective dimension at most $n$.
Hence we get that $M$ is of projective dimension at most $n$.\qed

\begin{lemma}\emph{\cite[Lemma 2.2]{JJ2}}\label{y3}
If $M\in\C$ has no direct summands in $\add \Sigma^nX$, then there exists an $(n+2)$-angle
$$X_0\xrightarrow{~~}X_1\xrightarrow{~~}\cdots\xrightarrow{~~}X_n\xrightarrow{~~}M\xrightarrow{~~}\Sigma^n X_0$$
in $\C$ with the following properties: Each $X_i$ is in $\add X$, and applying the functor ${\rm Hom}_\C(X,-)$ gives a
complex
$${\rm Hom}_\C(X,X_0)\xrightarrow{~}{\rm Hom}_\C(X,X_1)\xrightarrow{~}\cdots\xrightarrow{~}{\rm Hom}_\C(X,X_n)\xrightarrow{~}{\rm Hom}_\C(X,M)\xrightarrow{~}0$$
which is the start of the augmented minimal projective resolution of ${\rm Hom}_\C(X,M)$.
\end{lemma}

\section{Proof of the Main Theorem}
Now we begin to prove the main result of this article.

\begin{theorem}\label{main}
Let $\C$ be an $(n+2)$-angulated category with an Opperman-Thomas cluster tilting object $X$, and $\Gamma$ be the endomorphism algebra of $X$. Let $M$ be an indecomposable object in $\C$ which does not belong to $\add \Sigma^nX$. Then the $\Gamma$-module ${\rm Hom}_{\C}(X,M)$ is of infinite
projective dimension if and only if the factorization ideal $I_M$ is non-zero.
\end{theorem}

\proof Since $M$ is an indecomposable object in $\C$ which does not belong to $\add \Sigma^nX$,
then there exists an $(n+2)$-angle
\begin{equation}\label{t1}
\begin{array}{l}
X_0\xrightarrow{~\alpha_0~}X_1\xrightarrow{~\alpha_1~}X_2\xrightarrow{~~}\cdots\xrightarrow{~~}X_n\xrightarrow{~\theta~}M\xrightarrow{~\beta~}\Sigma^n X_0
\end{array}
\end{equation}
which satisfies the properties of Lemma \ref{y3}.

First, we assume $\Gamma$-module  ${\rm Hom}_{\C}(X,M)$  has infinite projective dimension. Applying the
functor $\Hom_{\C}(X,-)$ to the $(n+2)$-angle (\ref{t1}), we have the following exact sequence in $\mod\Gamma$:
$$(X,\Sigma^{-n}M)\xrightarrow{(X,~(-1)^n\Sigma^{-n}\beta)}(X,X_0)\xrightarrow{~}(X,X_1)\xrightarrow{~}\cdots\xrightarrow{~}(X,X_n)\xrightarrow{~}(X,M)\xrightarrow{~}0$$
where we omitted $\Hom_{\C}$ because of lack of space.
It follows that the morphism $$\Hom_{\C}(X,~(-1)^n\Sigma^{-n}\beta)\neq 0,$$ since otherwise the projective dimension of
${\rm Hom}_{\C}(X,M)$ would be at most $n$.
Choose a morphism $(-1)^n\Sigma^{-n}\gamma$ in $\Hom_\C(X,\Sigma^{-n}M)$ whose
image under $\Hom_{\C}(X,(-1)^n\Sigma^{-n}\beta)$ is non-zero, that is, the composition
$$X\xrightarrow{(-1)^n\Sigma^{-n}\gamma}\Sigma^{-n}M\xrightarrow{(-1)^n\Sigma^{-n}\beta}X$$
is non-zero. This yields the non-zero composition
$$\Sigma^nX\xrightarrow{~\gamma~}M\xrightarrow{~\beta~}\Sigma^nX.$$
Hence there exists a non-zero element $\beta\gamma$ in the factorization ideal $I_M$.
\medskip

Conversely, if the factorization ideal $I_M$ is non-zero, we prove that the
$\Gamma$-module ${\rm Hom}_{\C}(X,M)$ is not of projective dimension at most $n$. Otherwise, applying the
functor $\Hom_{\C}(X,-)$ to the $(n+2)$-angle (\ref{t1}), we have the following exact sequence in $\mod\Gamma$:
$$0\xrightarrow{~}\Hom_{\C}(X,X_0)\xrightarrow{\Hom_{\C}(X,~\alpha_0)}\Hom_{\C}(X,X_1)\xrightarrow{~}\cdots\xrightarrow{~}\Hom_{\C}(X,X_n)\xrightarrow{~}\Hom_{\C}(X,M)\xrightarrow{~}0.$$
Let $ba\colon\Sigma^nX\xrightarrow{a}M\xrightarrow{b}\Sigma^nX$ be any element in $I_M$.
Since $\Hom_{\C}(X,\Sigma^nX)=0$, we obtain that the morphism $b\theta\in\Hom_{\C}(X_n,\Sigma^nX)$ is zero.
Hence we have the following commutative diagram:
$$\xymatrix{0\ar[r]\ar[d]&\Sigma^nX\ar[d]^a\\
X_n\ar[r]^{\theta}\ar[d]&M\ar[d]^b\\
0\ar[r]&\Sigma^nX
}$$
By (N3), there exist morphisms $c\colon X\to X_0$
and $d\colon X_0\to X$ which make the
following diagram commutative:
$$\xymatrix{X\ar[r]\ar@{-->}[d]^c&0\ar[r]\ar@{-->}[d]&0\ar[r]\ar@{-->}[d]&\cdots\ar[r]&0\ar[r]\ar@{-->}[d]&0\ar[r]\ar[d]&\Sigma^nX\ar[d]^a\ar[r]^{(-1)^n}&\Sigma^nX\ar@{-->}[d]^{\Sigma^nc}\\
X_0\ar[r]^{\alpha_0}\ar@{-->}[d]^d&X_1\ar[r]^{\alpha_1}\ar@{-->}[d]&X_2\ar[r]\ar@{-->}[d]&\cdots\ar[r]&X_{n-1}\ar[r]\ar@{-->}[d]&X_n\ar[r]^{\theta}\ar[d]&M\ar[d]^b\ar[r]^{\beta}&\Sigma^nX_0\ar@{-->}[d]^{\Sigma^nd}\\
X\ar[r]&0\ar[r]&0\ar[r]&\cdots\ar[r]&0\ar[r]&0\ar[r]&\Sigma^nX\ar[r]^{(-1)^n}&\Sigma^nX}$$
Thus $\alpha_0c=0$ and then $\Hom_{\C}(X,c)=0$ since $\Hom_{\C}(X,\alpha_0)$ is a monomorphism.
Since the functor $\Hom_{\C}(X,-)$ is faithful, we have $c=0$ implies $\Sigma^nc=0$. Therefore $ba=\Sigma^nd\circ\Sigma^nc=0$.
 It follows that $I_M =0$, which is a contradiction to our
assumption.
By Corollary \ref{cor}, we get that the $\Gamma$-module ${\rm Hom}_{\C}(X,M)$  has infinite projective dimension.  
This finishes the proof. \qed
\medskip

As a special case of Theorem \ref{main} when $n=1$, we have the following.
\begin{corollary}\emph{\cite[Theorem 1.1]{BBT}}
Let $\C$ be a triangulated category with a cluster tilting object $X$, and $\Gamma$ be the endomorphism algebra of $X$. Let $M$ be an indecomposable object in $\C$ which does not belong to $\add \Sigma X$. Then the $\Gamma$-module ${\rm Hom}_{\C}(X,M)$ is of infinite
projective dimension if and only if the factorization ideal $I_M$ is non-zero.
\end{corollary}

\textbf{Panyue Zhou and Xingjia Zhou}\\
College of Mathematics, Hunan Institute of Science and Technology, Yueyang, Hunan, 414006, People's Republic of China.\\
E-mail: \textsf{panyuezhou@163.com and 1149077341@qq.com}

\end{document}